# DEFINING PROBABILITY DENSITY FOR A DISTRIBUTION OF RANDOM FUNCTIONS[1]

By Aurore Delaigle and Peter Hall

*University of Melbourne, University of Bristol, University of Melbourne and University of California at Davis*

The notion of probability density for a random function is not as straightforward as in finite-dimensional cases. While a probability density function generally does not exist for functional data, we show that it is possible to develop the notion of density when functional data are considered in the space determined by the eigenfunctions of principal component analysis. This leads to a transparent and meaningful surrogate for density defined in terms of the average value of the logarithms of the densities of the distributions of principal components for a given dimension. This density approximation is estimable readily from data. It accurately represents, in a monotone way, key features of small-ball approximations to density. Our results on estimators of the densities of principal component scores are also of independent interest; they reveal interesting shape differences that have not previously been considered. The statistical implications of these results and properties are identified and discussed, and practical ramifications are illustrated in numerical work.

**1. Introduction.** The concept of probability density for a random function is becoming increasingly important in functional data analysis. For example, it underpins discussion of the mode of the distribution of a random function, addressed in particular by Gasser, Hall and Presnell (1998), Hall and Heckman (2002) and Dabo-Niang, Ferraty and Vieu (2004a, 2004b, 2006). Nonparametric or structure-free methods for curve estimation from functional data involve the concept of density, not least because they generally are based on estimators of Nadaraya–Watson type which require di-

Received May 2009; revised August 2009.

[1]Supported by a grant from the Australian Research Council and a grant from the National Science Foundation.

vision by an estimator of a small-ball probability. See, for example, Ferraty, Goïa and Vieu (2002a, 2002b, 2007a, 2007b), Ferraty and Vieu (2002, 2003, 2004, 2006a, 2006b) and Niang (2002). There is of course a more general and very large methodology for functional data analysis, accessible via the monographs by Ramsay and Silverman (2002, 2005).

In this paper we take up directly the notions of the density and mode of the distribution of a random function. We argue that, while a density function is generally not well defined in this context, it is possible to define a meaningful concept of density for a specific scale or resolution level which is intrinsically linked to a particular dimension in a principal component representation. The challenge is to determine the dimension. We give an argument which leads directly from scale to dimension, through a simple approximation to a small-ball probability at a given scale.

The density approximation suggests a simple and appealing definition of mode and leads directly to an empirical approximation to density for a given dimension. Likewise, the approximation also enables two functions to be compared on the basis of their "relative likelihoods," that is, the heights of the density at the respective functions; although we shall not explore that feature in the present paper. We develop theoretical arguments describing both the approximation and the estimation of principal component score densities, and we give numerical illustrations of our conclusions.

Our empirical methods involve estimating the densities of principal component scores using approximations to those scores based on estimators of eigenvalues and eigenfunctions. This problem is itself of intrinsic interest, not least because principal component score densities reveal interesting shape differences. Our theoretical results describe properties of density estimators in this context.

The problem of determining the intrinsic dimension of the distribution of a random function for given scale is related to that of estimating the effective dimension of a sample of $p$-vectors when the sample size, $n$, is much less than $p$. Indeed, the connection between very high-dimensional data problems, and problems involving functional data, is drawn explicitly by Leng and Müller's (2006) "stringing" method which permits a random function to be computed from a long data vector. Leng and Müller suggest that the effective dimension of the transformed data be computed using principal component methods which also underpin our analysis.

This paper is organised as follows. In Section 2 we define a notion of density (the "log-density") that can be used in the functional data context, and from there we define a modal function which can be used to measure central tendency. The log-density depends on the densities of the principal component scores, and in Section 3 we show how to estimate these densities and study theoretical properties of these estimators. In Section 4 we provide theoretical arguments that justify the use of log-density in the functional data context. In Section 5 we use our estimators of the surrogate density, of the densities of the principal components scores and of the modal curve, to analyse an Australian rainfall dataset. We also illustrate the methods on some simulated data. The proofs of the main results are gathered in Section 6.

## 2. Main results and their implications.

2.1. *Decomposition into principal components.* As a prelude to summarising our main results we briefly revise important properties of functional principal component decomposition. See Besse and Ramsay (1986), Ramsay and Dalzell (1991), Rice and Silverman (1991) and Silverman (1995, 1996) for early work on this topic. Let $X$ be a random function supported on a compact interval $\mathcal{I}$. If the covariance function of $X$ is positive definite, it admits the following spectral decomposition:

$$(2.1) \qquad K(s,t) \equiv \operatorname{cov}\{X(s), X(t)\} = \sum_{j=1}^{\infty} \theta_j \psi_j(s) \psi_j(t),$$

where the expansion converges in $L_2$ on $\mathcal{I}^2$, and $\theta_1 \geq \theta_2 \geq \cdots$ are the eigenvalues, with respective orthonormal eigenvectors $\psi_j$, of the linear operator with kernel $K$.

The functions $\psi_1, \psi_2, \ldots$ form a basis for the space of all square-integrable functions on $\mathcal{I}$, and, in particular we can write, for $X$ and any square-integrable function $x$ on $\mathcal{I}$,

$$X = \sum_{j=1}^{\infty} \theta_j^{1/2} X_j \psi_j, \qquad x = \sum_{j=1}^{\infty} \theta_j^{1/2} x_j \psi_j,$$

where the quantities $X_j = \theta_j^{-1/2} \int_{\mathcal{I}} X \psi_j$ and $x_j = \theta_j^{-1/2} \int_{\mathcal{I}} x \psi_j$ are the principal component scores (sometimes referred to as the principal components) corresponding to functions $X$ and $x$. If $E(X) = 0$ then the above decomposition of $X$ is termed the Karhunen–Loève expansion, or generalised Fourier expansion. The $X_j$'s are always uncorrelated (this follows from orthogonality of the $\psi_j$'s), and we shall assume that they are independent. This is exactly correct if $X$ is a Gaussian process, and it is almost always assumed to be the case in empirical or numerical work. In such cases, as here, independence is often interpreted pragmatically; it captures the main features of a population, allows relatively penetrating theoretical analysis and motivates simple, practical methodology such as the estimators explored in Section 5 below. Particularly in the infinite-dimensional setting of functional data analysis, it seems impossible to use effectively general models for random variables that are uncorrelated but not independent. Such an approach leads to cumbersome methods and does not seem to allow useful insight into theoretical properties.

2.2. *Log-density function.* Many descriptive and predictive functional data analyses lean heavily on properties of the space of principal component scores. For example, it is common to describe properties of a sample of curves in terms of the shapes of the eigenfunctions corresponding to the largest eigenvalues obtained by functional principal component analysis. Similarly, it seems natural to define density for functional data in terms of the densities $f_j$ of the principal component scores $X_j$, for example, via the product of the densities $f_j$ corresponding to the largest eigenvalues. The notion of log-density, which we shall define below, fills this role.

For $h > 0$, let $p(x \mid h) = P(\|X - x\| \le h)$ where $\|X - x\|$ denotes the $L_2$ distance between $X$ and $x$. We shall show in Section 4 that

$$(2.2) \qquad \log p(x \mid h) = C_1(r, \theta) + \sum_{j=1}^{r} \log f_j(x_j) + o(r),$$

where $r = r(h)$ diverges to infinity as $h$ decreases to zero, $f_j$ is the density of the $j$th principal component score, $x_j$ is the version of that score for the function $x$ and both $r$ and the constant $C_1$ depend on $h$ and on the infinite eigenvalue sequence, $\theta$ say. (Neither $r$ nor $C_1$ depends on $x$ or on the distributions of the principal component scores.) The term $\sum_{j \le r} \log f_j(x_j)$ in (2.2) typically diverges at rate $r$ as $r$ is increased, and in particular it is generally not equal to $o(r)$.

Result (2.2) implies that, for appropriate eigenvalue sequences,

(2.3)    the log-density $\ell(x \mid r) = r^{-1} \sum_{j \le r} \log f_j(x_j)$ captures the variation, with $x$, of the logarithm of the small-ball probability $p(x \mid h)$, up to and including terms of order $r$ and in particular gives rise to a remainder of strictly smaller order than $r$.

[We have divided by $r$ only to ensure that $\ell(x \mid r)$ remains bounded as $r$ increases which makes it easier to discuss in theoretical terms. Of course, division by $r$ does not alter the main features of $\ell$.] More explicitly, we shall prove in Section 4 that

$$(2.4) \qquad p(x \mid h) = C_2(r, \theta) \exp\{r\ell(x \mid r) + o(r)\},$$

where $C_2(r, \theta) = (h\pi^{1/2})^r \Gamma(\frac{1}{2}r + 1)^{-1} \prod_{j \le r} \theta_j^{-1/2}$ does not depend on $x$. That is,

$$(2.5) \qquad p(x \mid h) = \frac{(h\pi^{1/2})^r}{\Gamma(r/2 + 1)} \left\{ \prod_{j=1}^{r} \theta_j^{-1/2} f_j(x_j) \right\} \exp\{o(r)\}.$$

One implication of the approximation at (2.4) is that it allows us to extract a function of $x$, namely the log-density $\ell(x \mid r)$, which captures the first-order effect that $x$ has on $p(x \mid h)$. In other words, the log-density describes the main differences in sizes of small-ball probabilities for different values of $x$. Moreover, up to terms that are negligible relative to those captured by the log-density $\log p(x \mid h)$, $\ell(x \mid r)$ is a monotone increasing function of $p(x \mid h)$. Therefore, while $\ell(x \mid r)$ cannot, in general, be employed to compare densities for different random function distributions, it can be used as the basis for comparing density at different points $x$ for the same random function distribution and for dimension $r$. Another implication of (2.4), which we shall derive in the Appendix, is that a probability density function for $X$ does not exist.

The principal component scores $X_j$ and $x_j$ are related to the squared $L_2$ distance between $X$ and $x$ by the formula

$$(2.6) \qquad \|X - x\|^2 = \sum_{j=1}^{\infty} \theta_j (X_j - x_j)^2.$$

Although the $x_j$'s are obviously linked to the eigenvalues $\theta_j$, in terms of the way these two quantities influence the distance of $X - x$ from zero [see (2.6)], it can be seen from (2.5) that the $\theta_j$'s and $x_j$'s are largely disconnected in terms of the way they influence the probability that $X - x$ is only a small distance from zero. In particular, we could arbitrarily permute the densities $f_1, \ldots, f_r$ without influencing anything other than the $o(r)$ term on the right-hand side of (2.5). Additionally, (2.4) makes it clear that the densities of the principal component scores are important only through their aggregate, defined by $\ell(x \mid r)$ in (2.3), and are of relatively little individual relevance.

2.3. *Defining the mode.* The log-density can be used to define a notion of central tendency in a population of curves. In the literature, central tendency is sometimes measured by the mean function. While this quantity is very easy to calculate, and it is close to its analogous definition in finite-dimensional settings, it is well known that it is generally unsatisfactory as a measure of "average" in the context of functional data since it tends to average out most of the fluctuations. For example, when applied to non-registered data or to data which cannot be perfectly aligned the averaging process often results in a mean curve that does not share typical properties (such as oscillations) of the population of curves.

As an alternative, we propose representing central tendency by the "modal function,"

$$(2.7) \qquad x_{\text{mode}} = \sum_{j=1}^{\infty} \theta_j^{1/2} m_j \psi_j,$$

which, for each $j$, has the $j$th principal component score $x_j$ equal to the mode $m_j$ of $f_j$. In a finite sample, $x_{\mathrm{mode}}$ can be estimated by

$$(2.8) \qquad \hat{x}_{\mathrm{mode}} = \sum_{j=1}^{T} \hat{\theta}_j^{1/2} \hat{m}_j \hat{\psi}_j,$$

where $\hat{\psi}_j$ and $\hat{\theta}_j$ are estimators of $\psi_j$, and $\theta_j$, $\hat{m}_j$ is the mode of $\hat{f}_j$, $\hat{f}_j$ an estimator of $f_j$, and $T$ is a truncation point which grows with the sample size.

Of course, in the case of functional data as well as multivariate data, we could also use the median to measure central tendency. In the context of a random function $X$, the median curve can be defined (analogously to the spatial median) to be the function $x$ that minimizes $E\|X - x\|$ (note that the theoretical mean minimizes $E\|X - x\|^2$). In practice, this median function can be estimated from the data by an iterative algorithm described in Gervini (2008). It can be shown through experimentation that the median is not as susceptible to the problem of averaging fluctuations as the mean but often more susceptible than the mode. For instance, in problems where the population consists of two or more well-separated sub-populations, the mean and the median can represent a function that is central in a strict mathematical sense but not representative of any function in any sub-population; whereas the modal function will often represent the most likely function in one of the sub-populations, and therefore will be less abstract and more interpretable than the mean or the mode. Nevertheless, in important cases (e.g., when $X$ is a Gaussian process) the mean, median and mode are identical.

## 3. Estimation of density of principal components, and log-density estimation.

3.1. *Empirical estimation of the density of principal component scores.* In this section we show how to estimate the densities $f_j$ of the principal component scores. This result will be used to provide an estimator of the log-density $\ell$, but it is of intrinsic interest, since having access to the densities $f_j$ can also be very useful for descriptive analysis of functional data. The $f_j$'s contain indeed valuable additional information about the structure of the population, compared to just the $\theta_j$'s and $\psi_j$'s.

Starting from independent data $X_{[1]}, \ldots, X_{[n]}$ on $X$, compute

$$(3.1) \quad \hat{K}(s,t) = \frac{1}{n} \sum_{i=1}^{n} \{X_{[i]}(s) - \bar{X}(s)\}\{X_{[i]}(t) - \bar{X}(t)\} = \sum_{j=1}^{\infty} \hat{\theta}_j \hat{\psi}_j(s) \hat{\psi}_j(t),$$

where $\bar{X} = n^{-1} \sum_i X_{[i]}$, the expansion in (3.1) is the empirical analogue of that at (2.1) and the terms are ordered so that $\hat{\theta}_1 \geq \hat{\theta}_2 \geq \cdots$. Thus we are centring the data at the sample mean rather than at the true mean which is of course unknown. We interpret $\hat{\theta}_j$ and $\hat{\psi}_j$ as estimators of the eigenvalues $\theta_j$ and eigenfunctions $\psi_j$, respectively. (We use square-bracketed subscripts so as not to confuse the $i$th data value $X_{[i]}$ with the $i$th principal component score, $X_i$, of $X$.) See, for example, Ramsay and Silverman [(2005), Chapters 8–10]. Then we calculate approximations $\hat{X}_{[ij]} = \hat{\theta}_j^{-1/2} \int_{\mathcal{T}} (X_{[i]} - \bar{X}) \hat{\psi}_j$ to the principal components $X_{[ij]} = \theta_j^{-1/2} \int_{\mathcal{T}} (X_{[i]} - EX_{[i]}) \psi_j$. We define too $\hat{x}_j = \hat{\theta}_j^{-1/2} \int_{\mathcal{T}} (x - \bar{X}) \hat{\psi}_j$, an estimator of $x_j = \theta_j^{-1/2} \int_{\mathcal{T}} (x - EX) \psi_j$.

An estimator $\hat{f}_j$ of the probability density function $f_j$ of $\theta_j^{-1/2}(X_j - EX_j)$ can be computed using the following standard kernel methods:

$$(3.2) \qquad \hat{f}_j(u) = \frac{1}{nh} \sum_{i=1}^n W\left(\frac{\hat{X}_{[ij]} - u}{h}\right),$$

where $h$ denotes a bandwidth and $W$ is a kernel function. For an introduction to kernel density estimation [see, for example, Silverman (1986) and Wand and Jones (1995)]. The value of $h$ could be chosen using standard methods for random data, reflecting the fact that $\hat{X}_{[ij]}$, for $1 \leq i \leq n$ is an approximation to the independent sequence $X_{[ij]}$, $1 \leq i \leq n$.

Provided the $j$th eigenvalue $\theta_j$ is not equal to either $\theta_{j-1}$ or $\theta_{j+1}$, the estimators $\hat{\theta}_j$ and $\hat{\psi}_j$ are root-$n$ consistent for $\theta_j$ and $\psi_j$ (modulo a change of sign of $\psi_j$), respectively, and so $\hat{x}_j = x_j + O_p(n^{-1/2})$. Note too that $\bar{X}$ cancels from the numerator inside the kernel in the definition of $\hat{f}_j(\hat{x}_j)$, and in fact,

$$(3.3) \qquad \hat{f}_j(\hat{x}_j) = \frac{1}{nh} \sum_{i=1}^n W\left\{\frac{\int_{\mathcal{T}} (X_{[i]} - x) \hat{\psi}_j}{h \hat{\theta}_j^{1/2}}\right\}.$$

In Section 3.2 below we show that this estimator is first-order equivalent to its "ideal" counterpart,

$$(3.4) \qquad \bar{f}_j(x_j) = \frac{1}{nh} \sum_{i=1}^n W\left\{\frac{\int_{\mathcal{T}} (X_{[i]} - x) \psi_j}{h \theta_j^{1/2}}\right\},$$

which we would use if we knew $\theta_j$ and $\psi_j$. Properties of $\bar{f}_j(x_j)$, as an estimator of $f_j(x_j)$, can be worked out using standard arguments. In particular, $\bar{f}_j(x_j)$ has variance and bias asymptotic to $wf_j(x_j)/(nh)$ and $\frac{1}{2}w_2 f_j''(x_j)h^2$, respectively, where $w = \int W^2$ and $w_2 = \int u^2 W(u)\, du$.

Our estimator of the log-density $\ell(x \mid r)$, in (2.3), is given by

$$(3.5) \qquad \hat{\ell}(\hat{x} \mid r) = \frac{1}{r} \sum_{j=1}^{r} \log \hat{f}_j(\hat{x}_j).$$

An attractive feature of $\hat{\ell}(\hat{x} \mid r)$ is the ease with which it can be computed for a range of values of $r$.

3.2. *Theoretical properties.* Here we show that the estimators at (3.3) and (3.4) are uniformly first-order equivalent. Since the variance and bias of $\bar{\bar{f}}_j(x_j)$ are generally of exact orders $h^2$ and $(nh)^{-1}$, respectively, then first-order equivalence is attained if $\hat{f}_j - \bar{\bar{f}}_j = o_p\{(nh)^{-1/2} + h^2\}$. Result (3.10) below is a strong form of this property.

The conditions we impose are the following:

(3.6)  for each $C > 0$ and some $\delta > 0$, $\sup_{t \in \mathcal{I}} E\{|X(t)|^C\} < \infty$ and $\sup_{s,t \in \mathcal{I}: \, s \neq t} E[\{|s - t|^{-\delta} |X(s) - X(t)|\}^C] < \infty$;

(3.7)  for each integer $r \geq 1$, $\theta_k^{-r} E\{\int_{\mathcal{I}} (X - EX)\psi_k\}^{2r}$ is bounded uniformly in $k$;

(3.8)  there are no ties among the $j + 1$ largest eigenvalues;

(3.9)  the density $f_j$ of the $j$th principal component score is bounded and has a bounded derivative; the kernel $W$ is a symmetric, compactly supported probability density with two bounded derivatives; for some $\delta > 0$, $h = h(n) = O(n^{-\delta})$ and $n^{1-\delta}h^3$ is bounded away from zero as $n \to \infty$.

Note that the assumptions on $h$ in (3.9) permit a bandwidth of conventional size, that is, $h \sim \text{const.}\, n^{-1/5}$ for any positive constant. The theorem remains valid if $W$ is the standard normal density, but more generally, infinitely supported kernels require assumptions about the rate at which their tails decrease. The use of infinitely supported kernels would alleviate difficulties that might arise when calculating $\log \hat{f}_j(u)$. However, the main features of the log-density $\ell$ (e.g., its modes) are identical to those of its exponentiated form. Therefore, in practice, to describe the main properties of a sample of curves it is not necessary to calculate logarithms, and we can simply work with the product of the estimated densities $\hat{f}_j$.

Given a square-integrable function $x$ defined on $\mathcal{I}$, put $\|x\|^2 = \int_{\mathcal{I}} x(t)^2 \, dt$. For each $c > 0$, let $\mathcal{S}(c)$ denote the set of $x$ such that $\|x\| \leq c$. Recall that $\hat{f}_j(\hat{x}_j)$ and $\bar{\bar{f}}_j(x_j)$ can be interpreted as functionals of $x$, and are defined by (3.3) and (3.4), respectively.

THEOREM 3.1. *If (3.6)–(3.9) hold then, for all $c > 0$,*

$$(3.10) \qquad \sup_{x \in \mathcal{S}(c)} |\hat{f}_j(\hat{x}_j) - \bar{f}_j(x_j)| = o\{(nh)^{-1/2}\}.$$

For the sake of brevity the proof of Theorem 3.1 is omitted. It can be found in a longer version of this paper [Delaigle and Hall (2008)].

3.3. *Joint density estimation.* Part of the simplicity of the log-density estimator, at (3.5), is that it involves only the marginal principal component score density estimators, $\hat{f}_j$, and not estimators of the joint densities of those scores. Of course, this is a consequence of the assumption that the scores are independent, but without that assumption the working statistician would be faced with not only a substantially more complicated density approximation, but the need to estimate joint densities. The latter problem is itself very challenging, unless sample size is large, since the accuracy of nonparametric density estimators decreases rapidly as dimension increases. Therefore the estimators that are produced under the assumption of independence enjoy a simplicity, that is, in many cases, a prerequisite for practical implementation.

## 4. Theoretical studies leading to results (2.3)–(2.5).

4.1. *Assumptions.* Let $X$ and $x$ be, respectively, a random and a fixed function on $\mathcal{I}$, and let $X_1, X_2, \ldots$ and $x_1, x_2, \ldots$ be their scores, defined in Section 2.1. For simplicity we assume that $E(X) = 0$, but if this condition does not hold then the mean of $X$ can be incorporated into $x_j$ by adding a term $E(X_j)$. For each $j$, let $f_j$ be the density of $X_j$, and note that, by definition of the scores, the $X_j$'s are uncorrelated and have mean zero and variance 1; in this work we assume that they are independent. See the last sentence of Section 2.1 for discussion of this condition.

For $j = 1, 2, \ldots$, let $W_j = X_j - x_j$ and let $g_j$ denote the probability density of $W_j$. Thus, $W_1, W_2, \ldots$ are independent random variables and, for all real $w$, $g_j(w) = f_j(w + x_j)$. For a given sequence of $x_j$'s we assume that

$$(4.1) \qquad \sup_{j \geq 1} E(W_j^2) < \infty$$

and that the sequence $\theta_1 \geq \theta_2 \geq \cdots$ of eigenvalues associated with the covariance of the function $X$ are positive numbers such that

$$(4.2) \qquad \sum_{j=1}^{\infty} \theta_j < \infty.$$

Note that, by (4.1) and (4.2), the series $\sum_j \theta_j W_j^2$ converges with probability 1.

Suppose too that each $g_j$ is differentiable at the origin, with $g_j(0) \neq 0$, and define $\rho_j = g_j'(0)/g_j(0)$. We shall assume that the densities $g_j$ admit Taylor expansions about the origin. In particular, we ask that, for each $\lambda > 0$, there exist a finite constant $A(\lambda) > 0$ such that

$$(4.3) \quad \sup_{j \geq 1} |\rho_j| < \infty, \qquad \sup_{j \geq 1} \sup_{|w_j| \leq \lambda} w_j^{-2} |g_j(w_j)g_j(0)^{-1} - (1 + \rho_j w_j)| \leq A(\lambda).$$

To understand the type of conditions this requires of the $f_j$'s and $x_j$'s, assume that the densities $f_j$ all have two bounded derivatives, and that for each $\lambda > 0$,

$$(4.4) \quad \inf_{j \geq 1} \inf_{|u| \leq \lambda} f_j(u) > 0, \qquad \sup_{|u| \leq \lambda} \{|f_j'(u)| + |f_j''(u)|\} < \infty.$$

For example, (4.4) holds if the principal components of $X$ are identically distributed with a density that has a bounded second derivative and does not vanish on the real line. Then (4.3) holds with $g_j(w) = f_j(w + x_j)$ for any bounded sequence of real numbers $x_j$. The case of an unbounded sequence $x_j$ can also be treated, but rather than the more general conditions imposed above, it requires assumptions and arguments that are related to specific density types. Therefore we shall not develop that case here.

4.2. *Approximation of the small ball probability.* Our first result, Theorem 4.1 below, will underpin our approximations to the value of

$$(4.5) \quad p(x \mid h) = p(h) = P\left(\sum_{j=1}^{\infty} \theta_j W_j^2 \leq h^2\right)$$

as $h \downarrow 0$. [To derive (4.5) we used (2.6).]

Given $h$ and $\lambda$ satisfying $0 < h \leq \lambda \theta_1^{1/2}$, we shall suppose that $r = r(h) \geq 1$ has been chosen such that

$$(4.6) \quad \theta_r^{-1} h^2 \leq \lambda^2.$$

Define $S = h^{-2} \sum_{j \geq r+1} \theta_j W_j^2$, and let $G = G(\cdot \mid h, r)$ denote the distribution function of $S$. Our first approximation to $p(h)$, at (4.5), is given by $q(h)$, described in the following theorem. The proof of the theorem is given in Section 6.

THEOREM 4.1. *Assume (4.1)–(4.3), and that $r$ is chosen so that (4.6) holds. Then*

$$(4.7) \quad p(h) = \exp\{\omega(h, \lambda)\lambda^2\} q(h),$$

*where*

$$(4.8) \quad q(h) = \frac{(h\pi^{1/2})^r}{\Gamma(r/2 + 1)} \left\{\prod_{j=1}^{r} \theta_j^{-1/2} f_j(x_j)\right\} \int_0^1 (1-t)^{r/2} \, dG(t),$$

$|\omega(h, \lambda)| \leq B(\lambda)$ *and the function* $B(\lambda) > 0$ *is nondecreasing in* $\lambda$ *and does not depend on* $h$ *or* $r$.

Next we apply Theorem 4.1 to develop more specific approximations to $p(h)$ for small $h$. Our results depend on the rate of convergence of the sequence $\theta_j$ to zero, and we consider two cases. We shall say that a sequence is "superexponential" if

$$(4.9) \qquad \theta_{k+1}/\theta_k \to 0 \qquad \text{as } k \to \infty,$$

or equivalently, if $\theta_k^{-1} \sum_{j \geq k+1} \theta_j \to 0$. More generally, we shall say that the sequence is "exponential" if

$$(4.10) \qquad \theta_k^{-1} \sum_{j=k+1}^{\infty} \theta_j \text{ is bounded} \qquad \text{as } k \to \infty.$$

When the eigenvalues converge to zero at a slower rate, nonparametric methods, where the notion of a functional-data density is typically employed, have much lower performance and so are less attractive and less likely to be used.

We also need to define the effective dimension, $r = r(h)$, for a given value of scale, $h$. In the superexponential setting, if for some $s$ the value of $h^2/\theta_s$ is "sufficiently close to 1" then we should take $r = s$, but otherwise we should take $r$ to be the unique integer for which $\theta_{r+1} < h^2 < \theta_r$. More specifically,

$$(4.11)$$

there exists a sequence of positive constants $c_1, c_2, \ldots$, depending on the eigenvalue sequence $\theta_1, \theta_2, \ldots$ and diverging to infinity, such that, if (for a given $h$) there exists $s \geq 1$ such that $|\log(h^2/\theta_s)| \leq c_s$, then we take $r = r(h)$ to be the infimum of such values; and if no such $s$ exists then we take $r$ to be the value for which $\theta_{r+1} < h^2 < \theta_r$.

In the case of an exponential sequence we define

$$(4.12) \qquad r = r(h, \lambda) = \arg\max\{j : \theta_j^{-1} h^2 \leq \lambda^2\},$$

and, for $j = 1, 2$, we let $\delta_j(s, \lambda)$ denote a quantity which satisfies

$$(4.13) \qquad \lim_{\lambda \to \infty} \limsup_{s \to \infty} \delta_j(s, \lambda) = 0.$$

Put $\Theta_j = -\frac{1}{2}\log\theta_j$ and $\phi_j = \log f_j(x_j)$. The proof of the next theorem is given in Section 6. The first part of (4.14) below is identical to (2.5).

THEOREM 4.2. *Assume* (4.1)–(4.3). *In the superexponential case, for* $r$ *as at* (4.11),

$$(4.14)$$
$$p(h) = \frac{(h\pi^{1/2})^r}{\Gamma(r/2+1)} \exp\{o(r)\} \left\{ \prod_{j=1}^{r} \theta_j^{-1/2} f_j(x_j) \right\}$$
$$= \exp\left[ \frac{1}{2} r\{\log(2\pi e h^2) - \log r + o(1)\} + \sum_{j=1}^{r} (\Theta_j + \phi_j) \right]$$

*as $h \to 0$, and when $\theta_j$ is exponential, for $r$ as at (4.12),*

$$p(h) = \frac{(h\pi^{1/2})^r}{\Gamma(r/2+1)} \exp\{r\delta_1(r,\lambda)\} \left\{ \prod_{j=1}^{r} \theta_j^{-1/2} f_j(x_j) \right\}$$

(4.15)

$$= \exp\left[ \frac{1}{2} r\{\log(2\pi e h^2) - \log r + \delta_2(r,\lambda)\} + \sum_{j=1}^{r} (\Theta_j + \phi_j) \right],$$

*where $\delta_1(r,\lambda)$ and $\delta_2(r,\lambda)$ satisfy (4.13).*

Theorem 4.2 shows that, for appropriate eigenvalue sequences, the approximations given at (2.4) and (2.5) hold. These approximations are appropriate when the eigenvalue sequence $\theta_j$ decreases to zero at an exponential rate which is the most important case from a practical viewpoint.

4.3. *Other implications of the theorems.* The integer $r = r(h)$ represents the dimension in which we make an approximation to the small-ball probability $p(x \mid h)$ at scale, or resolution level, $h$. In particular, (2.3) links the notion of density to dimension rather than, as is more commonly the case, to small-ball radius. In theoretical terms the connection is expressed through simple formulae such as (2.4) or (2.5). From an empirical viewpoint, (2.3) suggests that, rather than attempt to estimate small-ball probabilities for different values of $h$, so as to get a good idea of the way in which the notion of density changes as scale becomes finer, we can instead estimate the values of $\ell(x \mid r)$ for different values of $r$.

Note that $\ell(x \mid r)$ can be interpreted for increasing finite values of $r$, but $\ell(x \mid \infty)$, which we might define by taking the limit as $r \to \infty$ in (2.3), does not necessarily exist. In particular, we could change any finite number of the densities $f_j$ without altering the definition of density on an infinitesimal scale. Therefore, density on an infinitesimal scale, that is, as $h \to 0$ or as $r \to \infty$, is not identifiable unless we assume a model which asserts sufficiently close connections between early densities and principal component scores, and later ones.

An example where $\ell(x \mid \infty)$ is often well-defined arises when $x$ is taken to be the modal function $x_{\text{mode}}$ at (2.7). In that case, in order for the value of the log-density $\ell(x \mid r)$ to be well defined as $r \to \infty$, it is necessary only that $r^{-1} \sum_j \log f_j(m_j)$ converge. In particular, this condition is satisfied trivially if all the distributions of principal component scores are identical.

More generally, the value of $r$ can be interpreted as the dimension of the scale space when the unit of scale is $h$. The need to take scale, or resolution level, into account when discussing the density of a random function reflects the importance of scale in other settings. For example, Chaudhuri and Marron ([1999](#), [2000](#)) and Godtliebsen, Chaudhuri and Marron ([2002](#)) took a scale-space view of function estimation, noting that the viewpoint was already commonly used in areas such as imaging. There, the authors argue that one can learn different properties of the population at each scale where larger scales explain the overall structure of the population, and finer scales help understand the finer structure.

For a given scale $h$, the results of Theorem [4.1](#), and in particular ([4.11](#)) and ([4.12](#)), indicate that the effective dimension $r$ satisfies

$$(4.16) \qquad h^2 \approx \theta_r.$$

Property ([4.16](#)) implies that, if we are considering two distinct random function distributions for which the respective eigenvalue sequences decrease at different rates, then, for a sufficiently small value of scale, $h$, the corresponding dimension, $r$, is greater in the case of the random function with the less rapidly decreasing eigenvalue sequence. Of course, this makes intuitive sense.

One could determine the value of $r$ empirically by using relatively conventional methods for dimension estimation. See, for example, Horn ([1965](#)), Velicer ([1976](#)), Zwick and Velicer ([1986](#)), Peres-Neto, Jackson and Somers ([2005](#)) and Hall and Vial ([2006](#)). Alternatively we could simply estimate $\ell(x \mid r)$ for an increasing sequence of values of $r$, accessing in this way information about how density changes as we increase dimension and learning different properties of the population for each $r$. Theoretical properties of empirical principal components are discussed by, for example, Hall and Hosseini-Nasab ([2006](#), [2009](#)).

## 5. Numerical studies.

5.1. *Australian rainfall data.* We applied our density and mode estimation methods to an Australian rainfall dataset, available at **http://dss.ucar.edu/datasets/ds482.1** and analysed by Lavery, Kariko and Nicholls ([1992](#)). The data consist of daily rainfall measurements between January 1840 and December 1990, at each of 191 Australian weather stations. The function $X(t)$ represents the rainfall at time $t$ where $t$ denotes the period that has passed, in a given year, at the time of measurement. Rainfall at time $t$ was averaged over the years for which the station had been operating with the aid of a local polynomial smoother passed through discrete observations. One weather station (the 190th station) was removed from the collection of 191 since its rainfall pattern was very different from those for all other stations. Figure [1](#) depicts the yearly rainfall curves of the remaining 190 stations. On the left we show those stations (usually located in the north) which exhibit a "tropical" pattern, that is, those where most rain fell in mid to late summer; and on the right we show the stations (usually in the south) where the majority of rain came in cooler months.

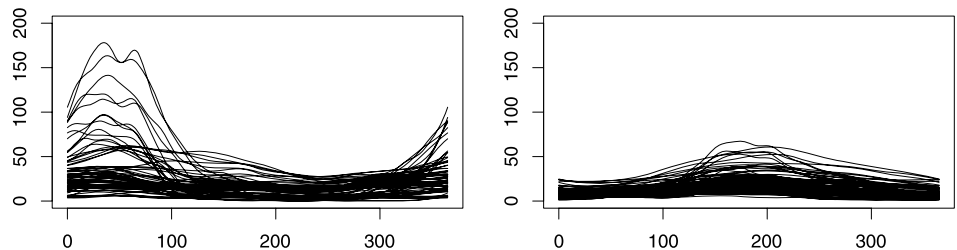

Fig. 1. *Australian rainfall data at weather stations which get the most rainfall during the summer months (left) or during the winter months (right).*

Table 1 indicates that the values of $\hat{\theta}_j$ decrease very quickly, and so it seems reasonable to use only the first few principal components to describe the data. In Figure 2 we plot the first four estimated principal component basis functions and the corresponding densities $\hat{f}_j$ of the standardised principal components scores. One interpretation of the first two principal component functions is that they capture, respectively, two key features of the data— peak rainfalls around days 35 and 215 and peak rainfalls around days 30 and 180. In each case, the first class of weather station generally corresponds to towns with a tropical or semi-tropical climate and the second to towns with a mid-latitude climate. Taken together, these two principal components capture the dichotomy between the two main latitude-determined climate zones in Australia together with the subtler effect of rainfall peaks that occur separately in either winter or summer, but where the peaks within either class can nevertheless differ by months.

We have chosen the signs of the first two principal component functions so that the first has its minimum in late winter whereas the minimum of the second is in late summer. However, this feature could easily be altered by a sign change; the signs of principal component functions are not determined. The densities of each of the first three principal component scores are skewed. Of course, the direction of skewness is tied to the sign of the principal component function which is arbitrary. The skewness is one aspect of the distinct non-Gaussian nature of the rainfall curves. The densities of higher-order principal component scores are less asymmetric; we plot only the fourth.

TABLE 1
*Proportion of variance explained by the first j principal components, for $j = 1, \ldots, 10$ in the Australian rainfall data example*

| j | 1 | 2 | 3 | 4 | 5 | 6 | 7 | 8 | 9 | 10 |
|---|---|---|---|---|---|---|---|---|---|----|
| | 0.7380 | 0.9510 | 0.9811 | 0.9915 | 0.9957 | 0.9974 | 0.9984 | 0.9989 | 0.9992 | 0.9995 |

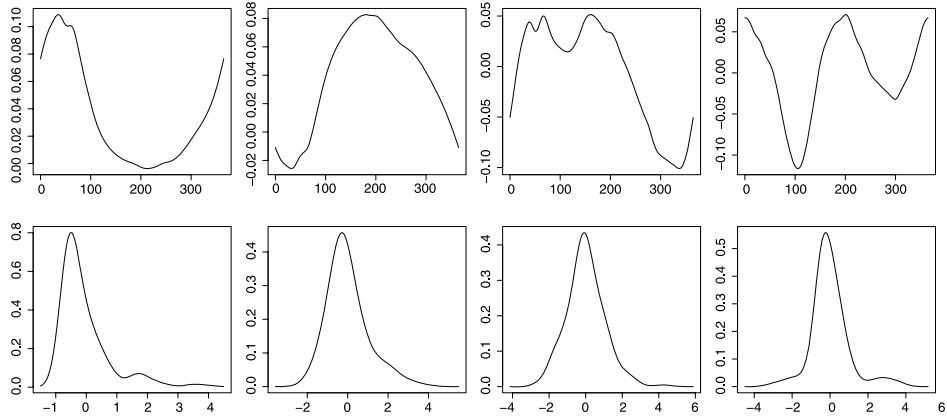

Fig. 2. *Plots of the first four principal component basis functions (row 1) and the corresponding estimates of densities of principal component scores, $\hat{f}_j$ (row 2). In each row, the graphs from left to right are for decreasing values of $\hat{\theta}_j$.*

Figure 3 shows the estimated mean and modal functions of the rainfall data where the modal function was calculated as at (2.8) and therefore depended on the number, $T$, of components used. Of course, since $\hat{x}_{\text{mode}}$ at (2.8) estimates the mode of the centered data, we have added the mean function to each modal curve. (The mean function is represented by the heavy, unbroken curve in Figure 3.) It is clear from that figure that changing $T$ from 1 to 5 alters the modal function significantly, but the effect of changing $T$ from 6 to higher values is almost indistinguishable by eye. The mean curve appears to be strongly influenced by the few stations that have high rainfall, but the modal curve is noticeably more robust. The figure also shows the median curve which we calculated using the algorithm described in Gervini (2008). The median lies between the mean and the mode, a property which is known to hold for many univariate distributions [see Haldane (1942) and Hall (1980)] but has not been studied previously for functional data. In this example the median has similarities with the modal curve, but it is still a bit high due to the influence of the tropical weather stations, especially in the summer months. These features make the median curve less appealing then the modal curve which looks more typical of curves for a majority of towns with mid-latitude climates.

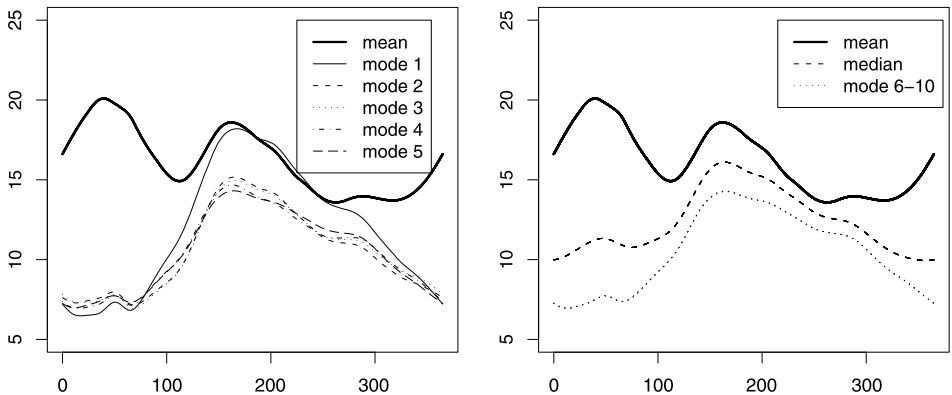

Fig. 3. *Curves representing the mean and mode, respectively, of Australian rainfall data, using $T = 1, 2, 3, 4$ or 5 (first panel) or $T = 6, 7, 8, 9$ or 10 (second panel). The annotation "mode $j$" indicates $\hat{x}_{\text{mode}}$ at (2.8) in the case $T = j$. On the right panel we also show the median curve.*

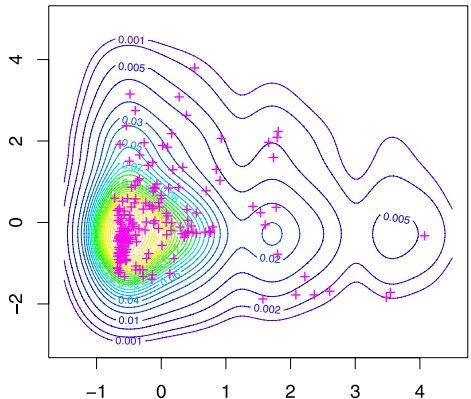

FIG. 4.   *Contour plot of the estimated surface* $2\exp\{\hat{\ell}(\hat{x} \mid 2)\} = \hat{f}_1(\hat{x}_1)\hat{f}_2(\hat{x}_2)$ *for the Australian rainfall data. The pink crosses represent the values of* $\hat{f}_1(\hat{X}_{[i1]})\hat{f}_2(\hat{X}_{[i2]})$ *for* $i = 1, \ldots, n$, *where* $X_{[ij]}$ *is the jth centered and scaled principal component score of the ith data curve.*

Figure 4 shows a contour plot of $2\exp\{\hat{\ell}(\hat{x} \mid 2)\} = \hat{f}_1(\hat{x}_1)\hat{f}_2(\hat{x}_2)$, that is, the estimated surrogate density for $r = 2$ together with the values $\hat{f}_1(\hat{X}_{[i1]})\hat{f}_2(\hat{X}_{[i2]})$, for $i = 1, \ldots, n$, represented by pink crosses. The colors range from blue for low density values to yellow for high density values.

We cannot visualise the estimated density $\hat{\ell}(\hat{x} \mid r)$ at higher resolution levels (i.e., for $r \geq 3$) by showing a surface curve. To see the effect that increasing $r$ has on $\hat{\ell}(\hat{x} \mid r)$, we calculated this density for $r = 1, \ldots, 10$ and for $x = X_{[1]}, \ldots, X_{[n]}$ (i.e., for each data curve), and then, for each $r$, classified the $n$ data curves into several groups according to the value of $\hat{\ell}(\hat{X}_{[i]} \mid r)$, using the same color code as Figure 4, that is, using colors ranging from blue for the lowest values of $\hat{\ell}(\cdot \mid r)$ to yellow for the largest values. We show in Figure 5 the groups of curves obtained for $r = 2$ and $r = 10$. We see that, overall, the curves of low (respectively, moderate or high) density for $r = 2$ correspond to the curves of low (respectively, moderate or high) density for $r = 10$. In other words, the density at resolution $r = 2$ already reflects the main features of the data. The blue curves roughly correspond to the stations for which rainfall varies the most over the year; these stations are very heterogeneous and thus have low density. At the other end of the spectrum, the yellow curves correspond to the stations with the flattest yearly rainfall; these stations are quite homogeneous and, logically, they have the highest density. The green curves correspond to a moderate rainfall change over the year and have moderate density values.

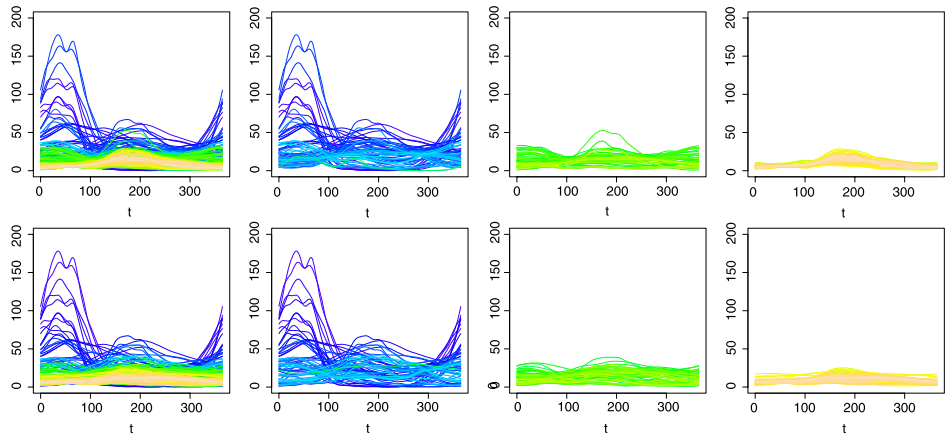

Fig. 5. *Plots of the 190 rain data curves. In each row, the furthest left graph shows all the curves and the other graphs show, from left to right, groups of curves for indices i that correspond to increasing values of $\hat{\ell}(X_{[i]} \mid r)$. The top row depicts results when $r = 2$ whereas the bottom row corresponds to $r = 10$. For $r = 2$ we use the color code corresponding to Figure 4, and for $r = 10$ we use a similar code corresponding to the estimator of $\hat{\ell}(X_{[i]} \mid 10)$.*

5.2. *Simulated examples.* As discussed at the end of Section 2.1, in the setting of functional data analysis it can be quite difficult to undertake meaningful inference without the simplifying condition of independence of the principal component scores $X_j$. In particular, without that assumption, to represent the joint density of $X_1, \ldots, X_r$ we would need to replace $\prod_{1 \le j \le r} f_j(x_j)$, a product of univariate functions, by a more complex $r$-variate function $f_{1,\ldots,r}(x_1, \ldots, x_r)$. However, the difficulty of estimating the latter increases rapidly with dimension. Therefore, unless samples sizes are particularly large, the quality of multivariate nonparametric estimators can be so poor that greater insight about the population is gained from estimators under the simplifying independence assumption. To illustrate this fact, we generated $B = 500$ samples of size $n = 100$ from the distribution of $X(t) = \sum_{1 \le j \le 10} \theta_j^{1/2} X_j \psi_j(t)$ where $t \in [0, 1]$, the $X_j$'s were uncorrelated dependent random variables generated according to $X_j = cTV_j$ where the $V_j$'s were independent and identically distributed, $T$ was a uniform $U[1, 2]$ variable, common to all $j$'s and $c = \{\mathrm{var}(TV_j)\}^{-1/2}$. We took $\psi_j(t) = \sqrt{2}\cos(\pi jt)$, and $V_j$ and $\theta_j$ were chosen from one of four models: (i) $V_j \sim \chi^2(8) - 8$, $\theta_j = j^{-3}$; (ii) $V_j \sim \chi^2(8) - 8$, $\theta_j = j^{-2}$; (iii) $V_j \sim \mathrm{N}(0, 1)$, $\theta_j = j^{-3}$; (iv) $V_j \sim \mathrm{N}(0, 1)$, $\theta_j = j^{-2}$ where "$\sim$" means "is distributed as."

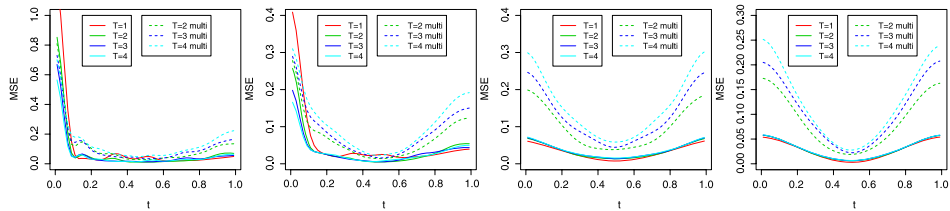

Fig. 6. *MSE of the estimators $\hat{x}_{\mathrm{mode}}$ and $\tilde{x}_{\mathrm{mode}}$ for $T = 1$, 2, 3 or 4. In the graphs we use* `multi` *to denote $\tilde{x}_{\mathrm{mode}}$. The graphs show, from left to right, the results of models* (i) *to* (iv).

In each case we calculated the estimator of the modal function, $x_{\mathrm{mode}}$. We compared $\hat{x}_{\mathrm{mode}}$ at (2.8) where each $f_j$ was estimated by a univariate kernel density estimator with plug-in bandwidth using the functions `kde` and `hpi` in Duong, Wand and Chacón's R package `ks` with the estimator $\tilde{x}_{\mathrm{mode}} = \sum_{1 \le j \le T} \hat{\theta}_j^{1/2} \tilde{m}_j \hat{\psi}_j$ where $(\tilde{m}_1, \ldots, \tilde{m}_T)$ was the mode of a $T$-variate kernel density estimator with plug-in bandwidth, calculated using the functions `kde` and `Hpi` in the R package `ks`.

In each case we calculated $\mathrm{MSE}(t) = B^{-1} \sum_{1 \le b < B} \{\hat{y}_b(t) - x_{\mathrm{mode}}(t)\}^2$ for $\hat{y}_b = \hat{x}_{b,\mathrm{mode}}$ and $\hat{y} = \tilde{x}_{b,\mathrm{mode}}$ where the index $b$ indicates that the estimator was calculated from the $b$th sample. In Figure 6 we present the MSEs of both estimators of $x_{\mathrm{mode}}$ for $T = 1$, 2, 3 and 4. The estimator $\hat{x}_{b,\mathrm{mode}}$ performed better than $\tilde{x}_{b,\mathrm{mode}}$ in all cases. Indeed, the quality of $\tilde{x}_{b,\mathrm{mode}}$ deteriorated very quickly as $T$ increased, due to greater bias and increased stochastic error; both these difficulties reflect the curse of dimensionality. For models (i) and (ii), the best results were obtained by $\hat{x}_{b,\mathrm{mode}}$ with $T = 3$ or 4, whereas for models (iii) and (iv), the best results were for $\hat{x}_{b,\mathrm{mode}}$ with $T = 1$. This reflects the fact that for models (i) and (ii), truncating the sum at $T < 10$ in the definition of $x_{\mathrm{mode}}$ introduced a systematic bias whereas in models (iii) and (iv), each $m_j = 0$ and thus truncating the sum at $T < 10$ did not produce any bias.

## 6. Technical arguments.

### 6.1. *Proof of Theorem 4.1.* Define $S_1 = \sum_{j \le r} \theta_j W_j^2$, $S = h^{-2} \sum_{j \ge r+1} \theta_j W_j^2$, $h_t = (1 - t)^{1/2} h$ and

$$(6.1) \quad p_r(h) = P(S_1 \le h^2) = \int_{\sum_{j \le r} \theta_j w_j^2 \le h^2} \left\{ \prod_{j=1}^r g_j(w_j) \right\} dw_1 \cdots dw_r.$$

Let $G = G(\cdot \mid h, r)$ denote the distribution function of $S$. In this notation,

$$(6.2) \quad p(h) = P(S_1 + h^2 S \leq h^2) = \int_0^1 P(S_1 \leq h_t^2) \, dG(t) = \int_0^1 p_r(h_t) \, dG(t).$$

Property (4.3) implies that

$$(6.3) \qquad \prod_{j=1}^r g_j(w_j) = \left\{ \prod_{j=1}^r g_j(0) \right\} \exp\left\{ \sum_{j=1}^r \rho_j w_j + O\left( \sum_{j=1}^r |w_j|^2 \right) \right\},$$

uniformly in $w_1, \ldots, w_r$ such that

$$(6.4) \qquad\qquad\qquad \sup_{1 \leq j \leq r} |w_j| \leq \lambda.$$

If (4.6) holds then so too does (6.4), provided $\sum_{j \leq r} \theta_j w_j^2 \leq h^2$. Therefore, by (6.1) and (6.3),

$$p_r(h_t) = h^r \left\{ \prod_{j=1}^r g_j(0) \right\} \int_{\sum_{j \leq r} \theta_j w_j^2 \leq 1 - t} \exp\left\{ h \sum_{j=1}^r \rho_j w_j + h^2 a(w) \right\} dw,$$

(6.5)

where $w = (w_1, \ldots, w_r)$, the function $a$ does not depend on $t$ and, for a constant $C_1 > 0$,

$$(6.6) \qquad\qquad\qquad |a(w)| \leq C_1 \sum_{j=1}^r |w_j|^2$$

uniformly in $w$ for which

$$(6.7) \qquad\qquad\qquad \sum_{j=1}^r \theta_j w_j^2 \leq 1.$$

The constant $C_1$ depends on $\lambda$, of which it is a nondecreasing function.

If (4.6) and (6.7) hold, then

$$(6.8) \qquad h^2 \sum_{j=1}^r w_j^2 \leq \lambda^2 \theta_r \sum_{j=1}^r w_j^2 \leq \lambda^2 \sum_{j=1}^r \theta_j w_j^2 \leq \lambda^2.$$

These properties, (6.5) and (6.6) imply that, for a constant $C_2 > 0$,

$$(6.9) \qquad\qquad\qquad h^2 |a(w)| \leq C_2 \lambda^2$$

uniformly in $w$ for which (6.7) holds. Here, $C_2 = C_2(\lambda)$ is a nondecreasing function of $\lambda$.

Let

$$(6.10) \qquad I = \int_{\sum_{j \leq r} \theta_j w_j^2 \leq 1 - t} dw = v_r (1 - t)^{r/2} \prod_{j=1}^r \theta_j^{-1/2},$$

where

$$(6.11) \qquad v_r = \frac{\pi^{r/2}}{\Gamma(r/2+1)}$$

denotes the content of the $r$-variate unit sphere. The second identity in (6.10) follows from the fact that, in view of the first identity, $I$ equals the content of the ellipsoid having the equation $\sum_{j \le r} \theta_j w_j^2 = 1 - t$.

Results (6.5) and (6.9) imply that

$$(6.12) \qquad p_r(h_t) = h^r \left\{ \prod_{j=1}^r g_j(0) \right\} \exp\{\omega_1(t \mid h, \lambda)\lambda^2\} J,$$

where, here and in (6.16) below, the function $\omega_j$ satisfies

$$(6.13) \qquad \sup_{0 \le t \le 1} |\omega_j(t \mid h, \lambda)| \le C_3$$

uniformly in $h$ and $\lambda$ such that (4.6) holds for some $r \ge 1$, the constant $C_3 > 0$ depends only on $C_2$ and the upper bound in (4.3) to $\sup_j |\rho_j|$, and

$$(6.14) \qquad J = \int_{\sum_{j \le r} \theta_j w_j^2 \le 1-t} \exp\left( h \sum_{j=1}^r \rho_j w_j \right) dw = I + \sum_{i=1}^\infty \frac{h^{2i}}{(2i)!} I_i$$

with

$$I_i = \int_{\sum_{j \le r} \theta_j w_j^2 \le 1-t} \left( \sum_{j=1}^r \rho_j w_j \right)^{2i} dw.$$

Odd-indexed terms have cancelled from (6.14) through symmetry.

When calculating $I_i$ using term-by-term expansion of the quantity within parentheses, only products of the form $(\rho_{j_1} w_{j_1}) \cdots (\rho_{j_{2i}} w_{j_{2i}})$, where each distinct index among $j_1, \ldots, j_{2i}$ appears an even number of times, make a nonzero contribution. Therefore,

$$(6.15) \qquad \begin{aligned} h^{2i} I_i &\le \int_{\sum_{j \le r} \theta_j w_j^2 \le 1-t} \left( h^2 \sum_{j=1}^r \rho_j^2 w_j^2 \right)^i dw \\ &\le (\rho\lambda)^{2i} \int_{\sum_{j \le r} \theta_j w_j^2 \le 1-t} dw = (\rho\lambda)^{2i} I, \end{aligned}$$

where $\rho = \sup_j |\rho_j|$ and, since (4.6) and (4.7) hold, we used the bound at (6.8). Using the bound (6.15) in (6.14) we deduce that $J = \exp\{\omega_2(t \mid h, \lambda)\lambda^2\} I$ where $\omega_2$ satisfies (6.13). This result and (6.12) imply that

$$(6.16) \qquad p_r(h_t) = h^r \left\{ \prod_{j=1}^r g_j(0) \right\} \exp\{\omega_3(t \mid h, \lambda)\lambda^2\} I.$$

The theorem follows on combining (6.2), (6.10), (6.11), (6.13) and (6.16).

6.2. *Proof of Theorem 4.2.* To derive (4.14) we treat two complementary cases which, if we consider convergence of $h$ to zero along subsequences, cover all instances: (a) there exists a sequence of integers $r = r(h)$ diverging to infinity such that, as $h$ converges along a subsequence, $h^2 \asymp \theta_r$ (that is, $h^2/\theta_r$ is bounded away from zero and infinity as $h \to 0$); and (b) along a subsequence, and for $r = r(h)$ diverging to infinity, $h^2/\theta_r \to 0$ and $h^2/\theta_{r+1} \to \infty$. In case (a) we take $\lambda^2$ to be an upper bound to $h^2/\theta_r$ and note that the superexponential condition (4.9) implies that $S \to 0$ in probability. From the latter property it follows that

$$(6.17) \qquad \frac{1}{r} \log\left\{ \int_0^1 (1-t)^{r/2} \, dG(t) \right\} \to 0.$$

In case (b) we choose $\lambda$ to be a function of $h$, decreasing to zero as $h \to 0$, in such a manner that $h^2 \theta_r^{-1} \lambda^{-2} \to 0$ and $h^2 \theta_{r+1}^{-1} \lambda^2 \to \infty$. The first of these properties ensures that (4.6) holds, and the second that $S \to 0$ in probability, so that, once again, (6.17) obtains. Hence, in either case, results (4.8) and (4.7), and Stirling's formula, imply (4.14). More generally, case (a) can be extended to that where $|\log(h^2/\theta_r)| \le c_r$, provided the sequence $c_1, c_2, \ldots$ diverges sufficiently slowly. This is the constant sequence used in the definition of $r(h)$ in (4.11).

To connect these results to the statement of (4.14) in Theorem 4.2, suppose that (4.14) fails in that context. Then we can find a sequence $h_1, h_2, \ldots$, decreasing to zero, such that the term written as $o(1)$ on the far right-hand side of (4.14) is actually bounded away from zero. Let $s = s(h)$ denote the integer for which $|\log(h^2/\theta_s)|$ is minimized. If, for all sufficiently large $k$, $|\log(h_k^2/\theta_{s(h_k)})| \le c_{s(h_k)}$, then $r(h_k) = s(h_k)$ [by the definition of $r(h)$ in (4.11)] and the result stated in the second-last sentence of the previous paragraph establishes (4.14). Hence, by passing to a sub-subsequence if necessary, we may assume that $|\log(h_k^2/\theta_{s(h_k)})| > c_{s(h_k)}$ for all sufficiently large $k$, in which case [again using the definition of $r(h)$] $\theta_{r(h_k)+1} < h_k^2 < \theta_{r(h_k)}$ for all large $k$, and both $h_k^2/\theta_{r(h_k)} \to 0$ and $h_k^2/\theta_{r(h_k)+1} \to \infty$. However, it then follows from case (b) in the previous paragraph that (4.14) holds. Therefore (4.14) must hold in the context of Theorem 4.2.

To prove (4.15), define $r$ as at (4.12) and note that $\theta_r^{-1} h^2 \le \lambda^2$ and $\theta_{r+1}^{-1} h^2 > \lambda^2$. The first of these properties ensures (4.6), and so permits us to apply Theorem 4.1, and the second guarantees that, with $C_1 = \sup_{j \ge 1} E(W_j^2)$ and $C_2$ denoting the upper bound to $\theta_k^{-1} \sum_{j \ge k+1} \theta_j$ in (4.10), we have

$$E(S) \le C_1 h^{-2} \sum_{j \ge r+1} \theta_j < C_1 \lambda^{-2} \theta_{r+1}^{-1} \sum_{j \ge r+1} \theta_j \le C_1 \lambda^{-2} (1 + C_2).$$

Therefore, given $\varepsilon_1 > 0$ we can choose $\lambda$ so large that $P(S > \varepsilon_1) < \varepsilon_1$, and hence, given $\varepsilon_2 > 0$ we can select $\lambda$ sufficiently large, but fixed, to ensure that, for all sufficiently small $h$,

$$(6.18) \qquad \left| \log\left\{ \int_0^1 (1-t)^{r/2}\, dG(t) \right\} \right| \le \varepsilon_2 r.$$

Results (4.8), (4.7), (6.18) and Stirling's formula imply (4.14).

## APPENDIX: NONEXISTENCE OF PROBABILITY DENSITY FUNCTION FOR FUNCTIONAL DATA

If $X$ is a random vector of finite length then we generally define the probability density, $f(x)$, of $X$ at the point $x$, as the limit as $h$ decreases to zero of the probability that $X$ lies in the ball of radius $h$ centered at $x$, divided by the Lebesgue measure of that ball. For example, in Euclidean space of dimension $r$,

$$(A.1) \qquad f(x) = \lim_{h \downarrow 0} (h^r v_r)^{-1} P(\|X - x\| \le h),$$

where $\|\cdot\|$ denotes Euclidean distance in $\mathbb{R}^r$, and $v_r$ represents the content of the $r$-dimensional unit sphere. It might be expected that a formula analogous to (A.1), with the divisor $h^r v_r$ replaced by a different function of $h$, would be appropriate for estimating the probability density of a random function $X$. However, in general it is not.

To appreciate why, let $X$ be a random function and $x$ a fixed function, and note that if there were to exist a function, $\alpha(h)$ say, such that the probability density

$$f(x) = \lim_{h \downarrow 0} \{\alpha(h)\}^{-1} P(\|X - x\| \le h)$$

were well defined, then for all $x$ we would have

$$\log f(x) = \lim_{h \downarrow 0} [-\log\{\alpha(h)\} + \log P(\|X - x\| \le h)]$$

and thus $\log P(\|X - x\| \le h) = C_1 + \log f(x) + o(1)$ where $C_1 = \log\{\alpha(h)\}$ does not depend on $x$, and $f(x)$ does not depend on $h$. However, (2.2) shows that this is not possible.

**Acknowledgments.** We thank two referees and an associate editor for their helpful comments which led to an improved version of the manuscript. The Australian weather data we used in the paper were assembled by the Australian Bureau of Meteorology. They are available from the Research Data Archive, maintained by the Computational and Information Systems Laboratory at the National Center for Atmospheric Research (NCAR). NCAR is sponsored by the National Science Foundation. Bob Dattore is acknowledged for providing the data.

Department of Mathematics
and Statistics
University of Melbourne
VIC, 3010
Australia
and
Department of Mathematics
University of Bristol
Bristol BS8 1TW
United Kingdom
E-mail: A.Delaigle@ms.unimelb.edu.au

Department of Mathematics
and Statistics
University of Melbourne
VIC, 3010
Australia
and
Department of Statistics
University of California
Davis, California 95616
USA
E-mail: halpstat@ms.unimelb.edu.au